\documentclass[11pt,reqno]{amsart}
\usepackage{graphicx,amscd}
\usepackage{amsfonts,amsmath,amssymb}
\newcommand{\rl}{{\mathbb{R}}}
\newcommand{\cx}{{\mathbb{C}}}

\newcommand{\G}{\Gamma}

\newcommand{\e}{\varepsilon}

\newcommand{\C}{\mathbb C}

\newcommand{\N}{\mathbb N}
\newcommand{\PP}{\mathbb P}

\newtheorem{theorem}{Theorem}
\newtheorem{lemma}[theorem]{Lemma}
\newtheorem{prop}[theorem]{Proposition}
\newtheorem{cor}[theorem]{Corollary}

\theoremstyle{definition}

\newtheorem{ex}[theorem]{Example}

\begin{document}
\title{On local hulls of Levi-flat hypersurfaces}
\author{Rasul Shafikov${}^{\alpha }$ and Alexandre Sukhov${}^{\beta, \gamma}$}
\begin{abstract}
We discuss local polynomial convexity of real analytic Levi-flat hypersurfaces in $\mathbb C^n$, $n>1$,
near singular points.
\end{abstract}

\maketitle

\let\thefootnote\relax

$\alpha$: Department of Mathematics, the University of Western Ontario, London, Ontario, N6A 5B7, Canada,
e-mail: shafikov@uwo.ca. The author is partially supported by the Natural Sciences and Engineering 
Research Council of Canada.

$\beta$: University of Lille, Departement de Math\'ematiques, 59655 Villeneuve d'Ascq, Cedex, France,
e-mail: sukhov@math.univ-lille1.fr. The author is partially supported by Labex CEMPI.

$\gamma$: Institut of Mathematics with Computing Centre - Subdivision of the Ufa Research Centre of Russian
Academy of Sciences, 450077, Chernyshevsky Str. 112, Ufa, Russia.


\section{Introduction}

Real analytic Levi-flat hypersurfaces naturally arise in the theory of holomorphic foliations, dynamical systems
 and some branches of complex geometry (such as minimal manifolds). From 
 this point of view, the most important problems are concerned with global properties of Levi-flat hypersurfaces. For example, 
 consider the existence problem: it is well-known, see~\cite{LN}, that in the complex projective space of dimension $> 2$
 there are no closed smooth real analytic Levi-flat hypersurfaces, but the question remains 
 open in dimension 2. The problem of existence disappears if one allows the hypersurface to have even mild singularities. 
 In general, real analytic Levi-flat hypersurfaces with singularities form a much more flexible class that can be constructed
 on many complex manifolds. This is one of the main 
reasons for the local study of singularities of Levi-flat hypersurfaces, a topic that has been recently developed by many authors . 
 This is also the principal motivation for our work. 

The present paper is devoted to the study of one of the crucial (local) complex analytic properties of Levi-flat hypersurfaces: 
polynomial or rational convexity near a singular point.
Local polynomial and rational convexity (or lack of it) is an important property of real submanifolds in complex Euclidean spaces 
that has many applications. It is well-known that any totally real submanifold is 
locally polynomially convex, while most CR submanifolds of positive CR-dimension are not. In this
paper we give a  characterization of local polynomial convexity for singular real analytic Levi-flat hypersurfaces. We refer the reader
to the next section for relevant definitions. 

In the classical local theory of real hypersurfaces (or, more generally, the Cauchy-Riemann manifolds) the properties of polynomially 
convex envelopes usually are determined from the intrinsic properties of the CR structure, such as nonvanishing of the Levi from 
(or its higher order analogs), or presence/absence of germs of complex analytic hypersurfaces, etc. In the Levi-flat  case these tools 
are not available: all Levi forms vanish identically, and the real hypersurface is foliated
by complex hypersurfaces (the Levi foliation). It turns out that local polynomial convexity is related to the geometric properties of the 
Levi foliation. This completely  new phenomenon is described in Theorem \ref{t.1}, which is our main result.

\begin{theorem}\label{t.1}
Let $\G\subset \mathbb C^n$, $n>1$, be an irreducible real analytic Levi-flat hypersurface, $0\in \G$. Then the following holds:
\begin{itemize}
\item[(i)] If $0$ is a regular point  of $\G$, or is an unbranched Segre nondegenerate singularity, then  $\G$ is locally polynomially convex at $0$.
\item[(ii)] There exists $\G \subset \cx^2$ with a two-branched Segre nondegenerate singularity at $0$ that is 
not locally rationally convex at $0$. Furthermore, there exists a neighbourhood basis $(U_k)_{k \in \N}$ of $0$ such that for every $k$ the polynomially convex hull of $\G  \cap U_k$
contains a full neighbourhood of $0$ in $\cx^2$.
\item[(iii)]  If $0$ is a Segre degenerate singularity, then $\G$ is not locally rationally convex at $0$, and
the rationally convex hull of  any compact neighbourhood $K\subset \G$ of $0$ in $\G$ contains a family of analytic  discs attached to $K$.
\end{itemize}
\end{theorem}

While a detailed presentation of the required background material is given in the next section, we would like to explain briefly the 
role of Serge (non)degeneracy and branching in the theorem. A singular point of a real analytic Levi-flat hypersurface
$\Gamma$ is called Segre degenerate, if the Segre variety of $\Gamma$ of this point coincides with $\C^n$. This case is settled by (iii). Otherwise the singularity is called Segre nondegenerate: its Segre variety  is a complex hypersurface. It is proved in \cite{SS} that the Levi foliation of $\Gamma$ extends to a full neighbourhood of a Segre nondegenerate singularity in $\C^n$ as a holomorphic web, that is, a (singular) holomorphic foliation with branching. It turns out that the local polynomial convexity near a Segre nondegenerate singularity depends on the degree of branching of this web. Part (i) says that when the web is not branched (i.e., it is a singular foliation), then we have local polynomial convexity (similarly to the nonsingular points). On the other and, in (ii) we present an example of a Segre nondegenerate 
singularity with a branched extension of the Levi foliation, which has a nontrivial hull.

In particular, it follows from Theorem~\ref{t.1} that for any point $p\in \G$, which is a unbranched Segre nondegenerate singularity, 
there exists a neighbourhood $K$ of $p$ in $\G$ such that there are no Riemann surfaces attached to $K$. We note that it is not true 
in general that a Levi-flat hypersurface (even everywhere smooth) does not admit ``large" holomorphic discs attached to it, as can be seen in the example of $\G=\{(z_1,z_2)\in \mathbb C^2 : |z_1|=1\}$. 
The following corollary characterizes locally polynomially convex points for a 
special class of Levi-flat hypersurfaces.

\begin{cor}\label{c.1}
Let $\Gamma\subset\mathbb C^n$, $n>1$, be a Levi-flat hypersurface such that its Levi foliation extends as a 
singular foliation to a neighbourhood of a singular point $p\in\G$. Then $\Gamma$ is locally 
polynomially convex at $p$ if and only if $p$ is a Segre nondegenerate singularity of $\Gamma$.
\end{cor}
This corollary can be applied to a wide class of Levi-flat hypersurfaces. For example, it was proved in \cite{SS} that the Levi foliation of  any algebraic (i.e. given by a polynomial defining function) Levi-flat hypersurface admits an extension as a holomorphic algebraic web.

We remark that  a singular point of $\G$ is Segre nondegenerate if and only if infinitely many leaves of the Levi foliation contain this point
in the closure, see~\cite{PSS1}. In dimension two such points are called {\it dicritical},  and so the above results can be reformulated 
using this terminology.

\medskip

A related problem is local holomorphic extension of CR-functions (see Section~\ref{s.convex} for the explicit connection).
By definition, CR-functions are those that satisfy the tangential Cauchy-Riemann equations. They are important in the 
context of boundary regularity of holomorphic functions and maps, and are ultimately related to 
singularities of holomorphic functions. Properties of CR-functions on CR-manifolds with singularities 
recently has attracted some attention, see, e.g.,
\cite{CS}, \cite{LNR} and references therein. While sufficient and necessary conditions are known for one-sided
holomorphic extension of all CR-functions from smooth real hypersurfaces, little is known in the case of hypersurfaces
with singularities. In Section~\ref{s.CR} we give a general criterion for extendability of CR functions for a certain class of 
(singular) hypersurface and give some examples of different behaviour of CR-functions on singular 
Levi-flat hypersurfaces.

\section{Background}

In this section we briefly review some relevant terminology and refer the reader to~\cite{BG}, \cite{CLN}, \cite{SS}, and~\cite{PSS1} 
for a detailed discussion of real analytic Levi-flat hypersurfaces. 

\subsection{Real analytic hypersurfaces}  Let $\Omega \subset \mathbb C^n$  be a domain. Consider 
a closed subset $\Gamma \subset \Omega$ such that locally, in a neighbourhood $U$ of  any point 
$q \in \Gamma$, it is given by 
\begin{eqnarray}
\label{hypeq}
\Gamma \cap U = \rho^{-1}( 0 ),
\end{eqnarray}
where $\rho: U \to \mathbb R$ is a real analytic function. We call $q \in \G$ a {\it regular} point, if $\G$ is a real analytic submanifold 
of dimension $2n-1$ in a neighbourhood of $q$, i.e., 
a smooth real analytic hypersurface near $q$. The union of all regular points forms the regular locus denoted by $\G^*$.  If $\G^*$ is not empty, then $\Gamma$ is called a real analytic hypersurface in $\Omega$. 
The set of regular points is open in $\Gamma$; its complement 
$\Gamma^{sng}:= \Gamma \setminus \Gamma^{*}$ is called the {\it singular locus} of $\Gamma$. 
The singular locus of $\G$ may contain points where $\G$ is a smooth manifold of dimension less than $2n-1$ (the so-called {\it stick} if $\G$ is irreducible). We emphasize that in this paper we ignore such points and simply call $\overline{\G^*}$ analytic hypersurface, even if this set is only semianalytic.  In other words, we identify $\G$ with $\overline{\G^*}$. 
From the point of view of foliation theory, this is natural  when one considers Levi-flat hypersurfaces, 
and this is our main motivation. With this convention, the regular locus
$\G^*$ is dense in $\G$.

We say that $\G$ is irreducible
if it cannot be represented as the union of two real analytic hypersurfaces. In fact, our considerations are purely local and we consider hypersurfaces irreducible as germs. We also assume that at the reference point, the defining function $\rho$ from (\ref{hypeq}) is minimal in the ideal of germs of real analytic functions vanishing on $\G$, see more details in \cite{PSS1}. 
However, for simplicity we will not use the terminology of germs in this paper.

Subanalytic sets are images of real analytic sets under proper real analytic maps. We refer the reader to~\cite{BM} 
for the theory of subanalytic sets.

\subsection{Levi-flat hypersurfaces and the Segre varieties}
 A hypersurface $\Gamma\subset \mathbb C^n_z$, $z=(z_1,\dots,z_n)$, is called Levi-flat if near every regular point it is 
 locally biholomorphically equivalent to a real hyperplane 
\begin{equation}\label{e.H}
H=\{ z \in \mathbb C^n: z_n + \overline{z}_n = 0 \} .
\end{equation}
This local equivalence induces a foliation on $\G^*$ by complex hypersurfaces, called the Levi foliation. 
Equivalently, $\G$ is Levi-flat if the restriction of the complex Hessian of  $\rho$ to the holomorphic tangent bundle  
of $\G^*$ (the Levi form) vanishes identically. 

For local analysis we may assume that $0\in \G$ and that $\varepsilon>0$ is so small
that the function $\rho$ in (\ref{hypeq}) admits the Taylor  expansion 
\begin{eqnarray}
\label{exp}
\rho(z,\overline z) = \sum_{IJ} c_{IJ}z^I \overline{z}^J, \ c_{IJ}\in\cx, \ \ I,J \in (\mathbb N \cup \{0\})^n ,
\end{eqnarray}
convergent in the ball $B(0,\e)$. Its complexification is defined by 
\begin{eqnarray}\label{compl}
\rho(z,\overline w) = \sum_{IJ} c_{IJ}z^I\overline{w}^J .
\end{eqnarray}
If $U$ is a sufficiently small neighbourhood of the origin, the series (\ref{compl}) converges for all $z,w \in U$. For $w \in U$ the complex analytic hypersurface, given by
\begin{eqnarray}
\label{Segre1}
Q_w = \{ z \in U : \rho(z,\overline w) = 0 \} ,
\end{eqnarray}
is called the {\it Segre variety} (associated with $\Gamma$) of the point $w$.  Segre varieties are defined invariantly 
with respect to the choice of the defining function $\rho(z)$ of $\G$, see~\cite{PSS1}. From the reality condition on $\rho$
it follows that 
$$
z \in Q_z  \,\, \mbox{ if and only if} \,\, z \in \Gamma,
$$
and
\begin{equation}\label{e.symm}
z \in Q_w \,\, \mbox{if and only if} \,\, w \in Q_z .
\end{equation}

\bigskip

  Let $q \in \Gamma^*$. Denote by ${\mathcal L}_q$ the unique leaf of the Levi foliation through $q$.
Then the leaf ${\mathcal L}_q$ is contained in the unique irreducible component of $Q_q$. 
In a small neighbourhood of $q$ this is also a unique complex hypersurface through $q$ which is contained in 
$\Gamma$. 

If $0\in\G$ is a singular point of $\G$, it may happen that in~\eqref{Segre1}, the function $\rho(z,0)$ vanishes identically, 
and so $Q_0=\mathbb C^n$. In this case we say that $0$ is a {\it Segre degenerate} point of $\G$. Segre degenerate
points form a complex analytic subset of dimension at most $n-2$. It follows from this that in $\mathbb C^2$
all Segre degenerate singularities are isolated. As was proved in~\cite{PSS1}, a point $p$ is Segre degenerate, if and only if 
infinitely many geometrically different leaves of the Levi foliation have $p$ in their closure. When $n=2$ this is equivalent
to $p$ being a dicritical singularity. 

\subsection{Singular webs and  foliations } In this subsection we recall the definition of a singular web, a detailed presentation is contained in \cite{SS}.  

We denote by $\PP T^*_n:=\PP T^*\C^n$ the projectivization of the cotangent bundle of $\C^n$ with the 
natural projection $\pi: \PP T^*_n \to \C^n$. A local trivialization of $\PP T^*_n$ is isomorphic to 
$U \times G(1,n)$, where $U\subset \cx^n$ is an open set and $G(1,n) \cong \mathbb CP^{n-1}$ is the Grassmannian 
space of linear complex one-dimensional subspaces in~$\cx^{n}$. The space $\PP T^*_n$ has the canonical
structure of a {\it contact manifold}, which can be described (using coordinates) as follows. Let $z = (z_1,...,z_n)$ be the
coordinates in $\C^n$ and $(\tilde p_1,..., \tilde p_n)$ be the fibre coordinates corresponding to the basis of differentials 
$dz_1, \dots, dz_n$.   We may view $[\tilde p_1, \dots, \tilde p_n]$ as homogeneous coordinates on $G(1, n)$.
Then, in the affine chart $\{\tilde p_n \ne 0\}$, with inhomogeneous coordinates $p_j = \tilde p_j / \tilde p_{n}$, 
$j=1,\dots, n-1$, the 1-form
\begin{equation}\label{e.cont}
\eta = d z_n + \sum_{j=1}^{n-1} p_j dz_j
\end{equation}
is a local contact form. Considering all affine charts $\{\tilde p_j \ne 0\}$ we obtain a global contact structure.

Let $U$ be a domain in $\C^n$. Consider a complex purely n-dimensional analytic subset $W$  in 
$\pi^{-1}(U) \subset \PP T^*_n$. Suppose that the following conditions hold:

\begin{itemize}
\item[(a)] the image under $\pi$ of every irreducible component of $W$ has dimension $n$;
\item[(b)] a generic fibre of $\pi$ intersects $W$ in $d$ regular (smooth) points and at every such point $q$ the differential 
$d \pi(q): T_qW \to \C^n$ is surjective;
\item[(c)] the restriction of the contact form $\eta$ on the regular part of $W$ is Frobenius integrable. So $\eta\vert_W = 0$ 
defines the foliation ${\mathcal F}_W$ of the regular part of $W$. (The leaves of the foliation ${\mathcal F}_W$ are called 
Legendrian submanifolds.)
\end{itemize}
Under these assumptions we define a {\it singular $d$-web} ${\mathcal W}$ in $U$ as a triple $(W,\pi,{\mathcal F}_W)$. 
A leaf of the web ${\mathcal W}$ is a component of the projection of a leaf of ${\mathcal F}_W$ into $U$. Note that at a 
generic point $z\in U$ a $d$-web $(W,\pi,{\mathcal F}_W)$ defines in $U$ near $z$ exactly $d$ families of smooth 
foliations. 

In dimension $2$ there is an immediate connection between singular webs and ODEs, which allows one to determine the value 
of the integer $d$. To describe that let $U \subset \cx^2$ be a domain, and  consider a holomorphic function $\Phi$  on 
$U \times \C$. It defines a {\it holomorphic ordinary differential 
equation} on $U \times \C$, 
\begin{eqnarray}\label{ME1}
\Phi(z_1,z_2,p) = 0
\end{eqnarray}
with $z = (z_1,z_2) \in U$ and $p= \frac{dz_2}{dz_1} \in \C$. This is an equation for the unknown function 
$z_2 = z_2(z_1)$.  
Any singular holomorphic $d$-web, $d \in \mathbb N$, can be defined in $U$ by equation (\ref{ME1}),  
 where $\Phi$ is of the form
\begin{eqnarray}\label{ME2}
\Phi(z,p) = \sum_{j=0}^d \Phi_j(z) p^{j}.
\end{eqnarray}
The graphs of solutions of (\ref{ME1}) are the {\it leaves} of  $\mathcal W$.  

\subsection{Extension of the Levi foliation and a meromorphic first integral}  
Let $\Gamma$ be a real analytic Levi-flat hypersurface in a domain $\Omega \subset \C^n$.
We say that a holomorphic $d$-web $\mathcal W$ in $\Omega$ is the {\it extension} of the Levi foliation of $\Gamma^*$
if every leaf of the Levi foliation  is a leaf of $\mathcal W$. In particular, if $d=1$, this defines the extension as a foliation,
singular in general. We assume that at least one leaf of every component of $\mathcal W$ agrees with a leaf of the Levi foliation;
under this condition the singular web extending the Levi foliation is unique.

The main result of \cite{SS} states that if $0$ is a Segre nondegenerate singularity of a real analytic Levi-flat hypersurface $\G$,
then the Levi foliation of $\G$ extends as a $d$-holomorphic web to a full neighbourhood of $0$ in $\cx^n$. The same conclusion 
holds if $\G$ is real algebraic, i.e., it is the zero-set of a real polynomial.
Thus, every singular point of $\G$ can be prescribed an integer $d$: we say that the origin is a {\it $d$-branched Segre nondegenerate singularity}. Section~\ref{s.BR} below contains an example of a Levi-flat hypersurface with a 2-branched singularity.
If $d=1$ we say that the origin is an {\it unbranched Segre nondegenerate singularity}, in this case the Levi foliation extends as a 
singular foliation to a neighbourhood of $0$ in the ambient space.

We also need a related notion of a {\it multiple-valued meromorphic first integral}.
Let $X$ and $Y$ be two complex manifolds and $\pi_X: X \times Y \to X$ and $\pi_Y :X \times Y \to Y$ be the natural 
projections. A $d$-valued {\it meromorphic correspondence} between $X$ and $Y$ is a complex analytic subset 
$Z \subset X \times Y$ such that the restriction $\pi_X\vert Z$ is a proper surjective generically $d$-to-1 map. 
Hence, $\pi_Y \circ \pi_X^{-1}$ is defined generically on $X$ (i.e., outside a proper complex analytic subset in $X$), 
and can be viewed as a $d$-valued map. In what follows we denote a meromorphic correspondence by a triple $(Z;X,Y)$ 
equipped with the canonical projections $\pi_X: Z \to X$ and $\pi_Y: Z \to Y$.

 A {\it multiple-valued meromorphic first integral} of a singular $d$-web $\mathcal W$ in $U$ is a $d$-valued meromorphic correspondence $(Z;U,\mathbb CP)$ with the following property: for a generic $c \in \mathbb CP$, 
 the set $R_c = \pi_U \circ \pi_{\mathbb CP}^{-1}(c)$ consists of a finite collection of complex hypersurfaces, 
and one of the irreducible components of $R_c$ agrees with some leaf of the Levi foliation.
It is proved in \cite{SS} that a $d$-web extending the Levi foliation near a Segre nondegenerate singularity always admits a 
meromorphic first integral. If the singularity is unbranched, then the extending web is a usual singular foliation and 
the first integral is a (single-valued) meromorphic function which is constant on every leaf.

We briefly recall this construction here to explicitly formulate the result needed in this paper. Let $\G$ be an irreducible Levi-flat
hypersurface with a Segre nondegenerate singular point $0\in \G$. Then there exists a complex line $A\subset \mathbb C^n$ such that
$Q_0\cap A = \{0\}$, $A \not\subset \G^{sng}$, and $A$ intersects every leaf of the Levi foliation near $0$.  For the unit disc 
$\mathbb D\subset \mathbb C$, let 
$$
\mathbb D \ni t \to w(\overline t) = (w_1(\overline t),\dots,w_n(\overline t)) \in \mathbb C^n
$$
be an anti $\mathbb C$-linear parametrization of the complex line $A$. 
The complexification of the defining function of $\G$ as in~\eqref{compl} defines a complex analytic hypersurface near the origin
in $\mathbb C^{2n}$, and this defines a complex analytic set
$$
Z = \{(z,t)\subset U \times \mathbb D : \rho(z,t)=\rho\left(z, \overline{w(\overline t)}\right)=0\} .
$$
By construction the natural projection $\pi: Z \to U$ is a proper map with discrete fibres. 
It is proved in~\cite{SS} that $\pi : Z \to U$ is 
the first integral of the singular web that extends the Levi foliation on $\G$. In particular, we obtain the following: 
{\it
If the origin is an unbranched  Segre nondegenerate singularity of an irreducible real analytic Levi-flat hypersurface $\G$, 
then there exists a 
holomorphic function $f: U \to V$, $V\subset \mathbb C$, which is constant along
the leaves of the Levi foliation on~$\G^*$.
}

\subsection{Polynomial, rational, and holomorphic convexity}\label{s.convex} 
Finally, we quickly recall the notion of polynomial and rational convexity.  
A compact $X\subset \mathbb C^n$ is called polynomially convex if the polynomially convex hull
$$
\widehat X = \{z\in \mathbb C^n : |P(z)| \le \sup_{w\in X} |P(w)|, P \text{ is any holomorphic polynomial} \}
$$
coincides with X. A compact $X$ is called rationally convex if it agrees with its rationally convex hull defined as the 
set of all points $z\in\mathbb C^n$ for which one cannot find a complex algebraic hypersurface that passes through $z$ and avoids $X$.
We say that $X$ is locally polynomially (resp. rationally) convex at a point $p \in X$, if there exists a neighbourhood basis of $p$ which consists of polynomially (resp. rationally) convex compacts. It is immediate from the definitions that if a compact 
$X\subset\mathbb C^n$ is not rationally convex, then it is not polynomially convex.

It is generally very difficult to determine whether a given compact is polynomially or rationally convex. One of the convenient tools is Oka's characterization of polynomial convexity.
The following formulation of Oka's principle can be found in Stout~\cite[Cor.~2.1.6]{St}: {\it Let $X\subset \mathbb C^n$ be a compact set, let $\Omega$ be an open set  that contains $\widehat X$􏰍, and let $z_0$ be a point of $\Omega$. 􏰅The point $z_0$ is not in $\widehat X$􏰍 if there exists a continuous family $\{V_t\}_{t \in [0,1)}$ of principal analytic hypersurfaces 
in $\Omega$􏰅 that diverges to infinity in $\Omega$􏰅 and that satisfies the conditions that $z_0\in V_0$ and $V_t\cap X =\varnothing$ for all $t\in [0, 1)$.} 

\smallskip

Let $X\subset \cx^n$ be compact and denote by $\mathcal O (X)$ the $\cx$-algebra of the (germs of) holomorphic functions on $X$.
We say that $X$  is  {\it holomorphically convex} if every nonzero $\cx$-algebra homomorphism 
$\chi : \mathcal O (X) \to \cx$ is of the form $\chi(f) = f(p)$ for some $p\in X$. Given a compact $X$ we say that the {\it holomorphic hull} of $X$ is the smallest holomorphically convex compact (with respect to inclusion) that contains $X$. If $M\subset \cx^n$ is a compact CR manifold, then the restriction to $M$ of any holomorphic function on $M$ is CR. Therefore, if every CR function on $M$ extends holomorphically to a fixed open set, then $M$ cannot be holomorphically convex, and therefore it has a nontrivial holomorphic hull. Clearly, this hull will be contained in the rationally convex hull of~$M$. We refer the reader to~\cite{St} for further discussion of holomorphically convex compacts.

\section{Proof of Theorem~\ref{t.1} and Corollary~\ref{c.1}}

\subsection{Proof of Theorem \ref{t.1}, (i)}
We first discuss the case of a regular point on $\G$. That a $C^1$-smooth real hypersurface is polynomially convex iff it is Levi-flat follows from Airapetyan~\cite{A}. For a smooth real analytic hypersurface $\G$ the proof of this is immediate. Indeed, if $p$ is a point on $\G$ where the Levi form of $\G$ does not vanish identically, then there exists a continuous family of holomorphic discs attached to $\G$ which contracts to $p$, see, e.g., Boggess~\cite{B} for details. This implies that $p$ is not polynomially (or even rationally) convex at $p$. The set of points where the Levi form has at least one nonzero eigenvalue is a dense opens subset of $\G$, and therefore, any neighbourhood $K$ of any point in $\G$ has Levi-nonflat points, 
and so $\G$ is not polynomially (and not rationally) convex at any point. The converse can be seen as follows: for any smooth point $p$ on a real analytic Levi-flat $\G$, there exists a neighbourhood of $p$ in $\G$ which is locally biholomorphically equivalent to the hyperplane 
$H$ given by~\eqref{e.H}. A convex compact in $H$ is polynomially convex, which can be easily seen from Oka's characterization of polynomial convexity, and therefore $\G$ is locally polynomially convex at $p$.

If $p\in \G$ is a singular point of $\G$ and $\G^*$ is not Levi flat, then any neighbourhood $K$ of $p$ in $\G$ contains points in $\G^*$ where the Levi form has at least one nonzero eigenvalue, and therefore $K$ is not polynomially (or rationally) convex by the argument above.

The remaining case is that $p=0$ is an unbranched Segre nondegenerate singularity of a Levi-flat $\G\subset \mathbb C^n$, $n>1$. Recall again that in this context we only deal
with $\G = \overline{\G^*}$ and ignore the points near which $\dim \G < 2n-1$. As discussed in the previous section, $\G$ admits a holomorphic first integral, i.e., there exist a open neighbourhood $U\in \mathbb C^n$ of $0$ and a holomorphic function $f: U \to \mathbb C$ that is constant on the leaves of the Levi foliation on $\G^*\cap U$. Let $Y=f^{-1}(f(\G))$. Then $f(\G)$ is a 1-dimensional subanalytic set and $Y$ is a subanalytic subset of $U$ of dimension $2n-1$ with $\G\subset Y$. Without loss of generality we may assume that $f(0)=0$. Since the set $f(\G)\subset \mathbb C$ is subanalytic of dimension 1,  there exists a neighbourhood $V$ of the origin in $\mathbb C$ such that $f(\G)\cap V$ consists of finitely many analytic arcs passing through the origin. In particular, the complement of $f(\G)$ in $V$ does not contain any connected components that are relatively compact in $V$.

Let $\e>0$ be so small that $B(2\e)\subset f^{-1}(V)$, and let $K=X \cap \overline{B(\e)}$.  Clearly, 
$\widehat K \subset B(2\e)$. To prove that $K=\widehat K$ we apply Oka's principle described above. 
Let $q\in B(2\e)\setminus Y$ be arbitrary. The set $f^{-1}(f(q))$ is a complex hypersurface in $U$. 
By construction, $f^{-1}(f(q)) \cap K =\varnothing$. To construct the required continuous family of complex hypersurfaces we simply choose a path $\gamma: [0,1] \to \mathbb C$ connecting the point $\gamma(0)=f(q)$ with a point $\gamma(1)$ outside $f(B(2\e))$ such that $\gamma$ avoids $f(\G)$. For each $t\in[0,1]$ the set $f^{-1}(\gamma(t))$ is then a complex hypersurface that avoids $K$ and 
$f^{-1}(\gamma(1))$ does not intersect $B(2\e)$. This shows that $q\notin \widehat K$. If now $q \in Y \setminus \G$, then the set $f^{-1}(f(q))$ contains several irreducible components at least one of which, say, $S_q$ passes through $q$. We claim that $S_q\cap \G = \varnothing$. Indeed, if $S_q$ intersects $\G$ then it intersects one of the leaves of the Levi foliation. From the positivity of the intersection index of complex varieties it follows that $S_q$ intersects all nearby leaves of the foliation on $\G$, but this is not possible because the leaves correspond to different level sets of $f$. As above we may find a path connecting $f(q)$ with a point outside $f(B(2\e))$ that avoids $f(\G)$, which again gives us the required continuous family of complex hypersurfaces. 
This shows that $K$ is polynomially convex, and so $\G$ is locally polynomially convex at the origin.

\subsection{Proof of (iii)}
Suppose now that $0\in X$ is a Segre degenerate singularity. It suffices to consider the case $n=2$, in which case
$0$ is an isolated dicritical singularity of $\G$.
Let $U$ be a neighbourhood of the origin where all Segre varieties are well-defined. We claim that there exists a point $w_0\in U \setminus \G$ such that $Q_{w_0} \cap \G = \{0\}$. Indeed, let 
$$\Sigma = \bigcup_{w\in \G\cap U \setminus\{0\}} Q_w .$$
Let $\G^c \subset \mathbb C^2_{z} \times \mathbb C^2_{w}$ be the complexification of $\G$, with the coordinate projections 
$\pi_z : \G^c \to \mathbb C^2_{z}$ and $\pi_w : \G^c \to \mathbb C^2_{w}$. Then 
$$
\Sigma = \pi_z (\pi_w^{-1}(\G)).
$$
Since the generic fibres of $\pi_w$ and $\pi_z$ have (complex) dimension $1$, it follows that 
$\dim_{\mathbb R} \pi_w^{-1}(\G) = 5$, and $\dim_{\mathbb R} \Sigma = 3$. Thus, there exists a point 
$w_0\notin \Sigma$. Then $Q_{w_0}$ passes through the origin (since $0$ is dicritical) but does not contain
any other points in $\G$. The latter can be seen as follows: if $w \in Q_{w_0} \cap \G$, $w\ne 0$, then $w_0\in Q_w$
by~\eqref{e.symm}, but this contradicts $w_0\notin \Sigma$. 

We now translate $Q_{w_0}$ slightly so that the translated hypersurface $S$ intersect smooth points on $\G$.
Then the set $S\cap \G$, which is a real analytic set (we may have to add points from the stick in $\G$ if it exists), 
contains a closed real curve. This curve bounds a domain in $S$ because it can be contracted to a point by shifting 
$S$ back to $Q_{w_0}$. This domain is a holomorphic disc which is attached to $\G$. This shows that no neighbourhood of 
the origin in $\G$ can be polynomially convex. 

Finally, no neighbourhood of the origin can be rationally convex for the following reason: let $z$ be a point on one of the discs attached to $\G$ through the above process. If $z$ is not in the rationally convex hull of a compact neighbourhood $K\subset \G$ of the origin, then there exists a complex algebraic hypersurface $R$ that passes through $z$ and avoids $K$. But $R$ intersects the variety attached to $\G$, and by the positivity of the intersection index for varieties it follows that $R$ either intersects the whole family of the varieties attached to $\G$, and hence passes through the origin, or $R$ intersects the boundary of some variety, i.e., intersects $K$. This contradiction shows that $K$ is not rationally convex.

\subsection{Proof of (ii)}\label{s.BR}
The following example, discovered by M.~Brunella~\cite{B1}, see also~ \cite{SS}, shows that, in general, 
the Levi foliation of a Levi-flat hypersurface admits extension to a neighbourhood 
of a singular point only as a web, not as a singular foliation. Consider the Levi-flat hypersurface
\begin{equation}\label{e.b1}
\G = \{z\in \mathbb C^2 : y_2^2=4(y_1^2+x_2)y_1^2\}.
\end{equation}
The singular locus of $\G$ is the set $\{y_1=y_2=0\}$. Its subset given by $\{y_1=y_2=0,\ x_2 <0\}$
is a stick, i.e., it does not belong to the closure of smooth points of $\Gamma$. After the complexification
we see that  $Q_0= \{z_2^2 +z_1^4 - 2 z_2 z_1^2=0\}$, and so the origin is a Segre nondegenerate singularity. Calculations in~\cite{SS} show that the singular web of $\Gamma$ at the origin in the form of representation~\eqref{ME1} can be given as
\begin{equation}\label{e.bb}
p^2 = 4z_2.
\end{equation}
This is precisely the 2-web that extends the Levi foliation of $\G^*$. Solving~\eqref{e.bb} shows that 
the first integral of $\G$ can be taken to be
$$f(z_1,z_2)=z_1 \pm \sqrt{z_2}.$$ 
Further, the closure of the smooth points of $\Gamma$ can be given by
$$
\{z\in \mathbb C^2 : {\rm Im\,}(z_1 \pm  \sqrt{z_2}) =0 \} =
\{{\rm Im\,}(z_1 +  \sqrt{z_2}) = 0\} \cup \{{\rm Im\,}(z_1 - \sqrt{z_2})=0\} .
$$

Consider now the holomorphic polynomial map $$F: (w_1,w_2) \mapsto (z_1,z_2) =  (w_1,w_2^2).$$ It is easy to see that 
$F^{-1}(\G) = \Pi_1 \cup \Pi_2$,  where $\Pi_1 = \{ {\rm Im\,}(w_1 + w_2) = 0 \}$ and 
$\Pi_2 = \{ {\rm Im\,}(w_1 - w_2) = 0 \}$. After a complex linear change of coordinates we may assume 
that the hyperplanes have the form 
$\Pi_j = \{ y_j = 0 \}$ and their intersection is $\rl^2$.

We employ a classical construction of complex discs. The domain $\cx^2 \setminus \rl^2$ is the union of 4 wedge type  domains
$\{ \tau_j y_j < 0, j= 1,2 \}$, where $\tau_j \in \{ -1, 1\}$. Consider, for example, the wedge $W = \{ y_j < 0, j= 1,2 \}$. It is contained in the strictly pseudoconvex domain
$$
\Omega = \{ y_1 + y_2 + y_1^2 + y_2^2 < 0 \} ,
$$ 
which is biholomorphic to the unit ball. Then there exists a complex curve touching the boundary $b\Omega$ of $\Omega$ from
outside exactly at the origin. Translating this curve to $\Omega$ in the direction of the inward normal to $b\Omega$ and considering the intersections of these complex curves with $W$, we obtain a family of complex discs filling $W$. This family contracts to the origin. Repeating this for other 3 wedges, we obtain the filling of a neighbourhood of $0$ in $\cx^2$ by a family of discs with boundaries contained in $\Pi_1 \cup \Pi_2$.

Now as in the proof of part (iii) it follows that $\G$ is not locally rationally 
convex at $0$.

\subsection{Proof of Corollary~\ref{c.1}}
 If $p$ is Segre degenerate, the result is immediate from Theorem~\ref{t.1}(iii). 
 Suppose that $p$ is Segre nondegenerate. It was proved in~\cite{CLN} and~\cite{B2} that if the Levi foliation admits,
locally near a singular point $p$, extension as a foliation, then near $p$ there exists a meromorphic first 
integral. On the other hand, by~\cite{SS}, there exists a $d$-valued first integral, which
by the uniqueness theorem (for example, for complex-analytic sets--graphs of the first integrals), must
be single-valued. Then by Theorem~\ref{t.1}(i), $\G$ is locally polynomially convex at $p$.

\section{Extension of CR-functions}\label{s.CR}

\subsection{Hans Lewy's extension principle}
The following general principle of holomorphic extension of CR functions from singular real analytic
hypersurfaces is inspired by H.~Lewy's method. It does not use the holomorphic approximation of 
CR functions.

Suppose that $\G$ is a real analytic hypersurface in $\cx^{n+1}$, $0\in \G$. Let $\Omega$ be a
one-sided neighbourhood of the origin, i.e., a connected component of $(\cx^{n+1}\setminus \G )\cap U$,
for some neighbourhood $U\ni 0$. Suppose that the coordinates
$(z,w)\in \cx^{n}_z\times \cx_w$ are chosen in such a way that the following conditions hold:

\begin{itemize}
\item[(i)] there exists a domain $D \subset \cx^n$ such that for every $z \in D$ the set
$$\gamma(z) = \{ w \in \cx: (z,w) \in \G \}$$
is a connected real analytic curve that bounds a simply connected domain 
$$
\Delta(z) = \{ w \in \cx : (z,w) \in \Omega \} .
$$
Further assume that $\gamma(z)$ is contained in a regular part of $\G$ for all $z\in D$.
\item[(ii)] the union $\Omega = \cup_{z \in D} \Delta(z)$ is a domain in $\cx^{n+1}$. Note that part of its boundary is contained in $\G$.
\end{itemize}
Consider a continuous function $f: \G \to \cx$ which is a CR function of class $C^1$ on the regular part of $\G$.
Denote by $E(f)$ the set of points $z \in \overline D$ such that one of the following conditions hold: 

\begin{itemize}
\item[(iii)] $\gamma(z)$ is a single point;
\item[(iv)] the function $\zeta \mapsto f(z,\zeta)$ extends holomorphically from $\gamma(z)$ to $\Delta(z)$.
\end{itemize}

\begin{prop}\label{t.2}
If $E(f)$ is the set of uniqueness  for functions holomorphic on $D$ and continuous on $\overline D$,
then $f$ extends holomorphically to $\Omega$.
\end{prop}

We note that in the case considered by H.~Lewy, $E(f)$ is an open arc on the boundary of $D$.
This can be seen, for example, in the case of the unit sphere $\G: \vert z \vert^2 + \vert w \vert^2 = 1$: the 
domain $D$ is the unit disc and $\gamma(z)$ shrink to a single point over its boundary. 

Our proof is based on the classical method of H.Lewy \cite{L}. It was applied by A.Tumanov \cite{T} for some classes of piece-wise real analytic Levi-flat hypersurfaces.

\begin{proof}
We prove that the function 
\begin{eqnarray}
\label{Cauchy}
F(z,w) = \frac{1}{2\pi i} \int_{\gamma(z)} \frac{f(z,\zeta) d\zeta}{\zeta - w}
\end{eqnarray}
defines the holomorphic extension of $f$ to $\Omega$. Clearly, $F$ is holomorphic in $w$. 
We also consider the moment functions
\begin{eqnarray}
\label{moment}
\mu_k(z) = \int_{\gamma(z)} \zeta^k f(z,\zeta)d\zeta, \,\,\, k = 0,1,2,\dots .
\end{eqnarray}

\begin{lemma}
The functions $F$ and $\mu_k$ are holomorphic in $z$ on $D$.
\end{lemma}

We prove the lemma for $F$; the argument for  $\mu_k$ is the same. By the Hartogs separate
analyticity theorem it suffices to prove the result on any slice $D_j = L_j \cap D$, $j=1,\dots, n$,
where $L_j\subset \cx^n$ is any complex line parallel to $\cx_{z_j}$.
Fix a simple closed smooth curve $\tau \subset D_j$ that bounds there a domain $D_\tau$, and consider in 
$\G \cap (L_j \times \cx_w)$ a $2$-torus given by
$$T = \{ (z,\zeta) \in \G: z \in \tau,\, \zeta \in \gamma(z) \}.$$
It bounds a solid torus $N \subset \G \cap (L_j \times \cx_w)$ of the form
$$N = \{ (z,\zeta) \in \G: z \in D_\tau, \zeta \in \gamma(z) \} .$$
Fix $w$ and set $G(z,\zeta) =  \frac{1}{2\pi i} \frac{f(z,\zeta)}{\zeta - w}$.
Using the assumption that $f$ is a CR function, we obtain by Stokes' formula
\begin{eqnarray*}
  \int_\tau F(z,w)d z_j = \int_T G(z,\zeta)d\zeta \wedge d z_j = \int_N dG(z,\zeta) \wedge d\zeta \wedge d z_j = 0 .
\end{eqnarray*}
By Morera's theorem we conclude that $F$ is holomorphic in $z_j$. This proves the lemma.

\medskip 

By assumption, for any $z\in E(f)$, there are two possiblities. First,  the function $f$ extends holomorphically to $\Delta(z)$. This implies that for all $k\ge 0$
the moment functions $\mu_k(z)$ vanish for these $z$. Another possiblility is that $\gamma(z$) degenerates to a single point. Then $\mu_k(\zeta)$ tends to $0$ as $\zeta \in D$ tends to $z$. Since $E(f)$ is the set of uniqueness we conclude $\mu_k(z)\equiv 0$ for
all $k$, which implies that $f$ extends holomorphically to $\Delta(z)$ for all $z\in D$.
This proves that $F$ gives a desired holomorphic extension of $f$.
\end{proof}

%

\subsection{Examples}
In this section we give examples of Levi-flat hypersurfaces with different extension property of CR-functions
defined on them.

\begin{ex}
Consider the Levi-flat hypersurface $\G$ given by
$$\G = \{ (z,w) \in \cx^2: \vert z \vert^2 - \vert w \vert^2 = 0 \}.$$
This is a real cone which divides the ambient space into two connected components. Note that the complex lines 
$\{z=0\}$ and $\{w=0\}$ intersect $\G$ at the origin only and belong to two different components of the 
complement of $\G$ (except $0$). In the terminology of~\cite{CS}, $\G$ is said to have two-sided support. 
It was proved in~\cite{CS} that for singular hypersurfaces with two-sided support there always exist CR-functions 
that are $C^m$-smooth (for any $m\in\mathbb N$) that do not extend holomorphically to either side of the hypersurface. 
For example, for the given $\G$ one can take 
$$
f(z,w) = \frac{z^3}{w}+\frac{w^3}{z}.
$$
This $f$ is differentiable on $\G$, CR on the regular part of $\G$ and clearly does not extend holomorphically 
to either side of $\G$.

Now let
$$\Omega = \{ (z,w) \in \cx^2: \vert z \vert > \vert w \vert \}$$
be one component of $\cx^2 \setminus \G$. The $z$-projection sends $\Omega$ to $\cx^* = \cx_z \setminus \{ 0 \}$. 
In the notation of Proposition~\ref{t.2}, for every 
$z \in \cx^*$ we have 
$$\Delta(z) = \{ w \in \cx : (z,w) \in \Omega \} ,$$
and
$$\gamma(z) = \{ w: \vert w \vert = \vert z \vert \}.$$
If there exists a set of uniqueness $E(f)\subset \cx^*$ such that a smooth CR-function $f$ on $\G\setminus\{0\}$ 
extends holomorphically to $\Delta(z)$ for all $z\in E(f)$, then by Proposition~\ref{t.2}, $f$ extends holomorphically
to~$\Omega$. In this example, one may take $E(f)$ to be, for instance, a sequence of points in $\cx^*$ that
has an accumulation point (in $\cx^*$). Also note that in view of the above discussion the requirement on the existence of $E(f)$
cannot be dropped in general. 
\end{ex}

\begin{ex} In this example we continue the exploration of the Levi-flat hypersurface given in  Section~\ref{s.BR}. Let
\begin{equation}\label{e.b1}
\G = \{z\in \mathbb C^2 : y_2^2=4(y_1^2+x_2)y_1^2\}.
\end{equation}
We recall that the closure of the regular points of $\G$ can be given as 
$$
\G^*=\{z\in \mathbb C^2 : {\rm Im\,}(z_1 \pm  \sqrt{z_2}) =0 \} =
\{{\rm Im\,}(z_1 +  \sqrt{z_2}) = 0\} \cup \{{\rm Im\,}(z_1 - \sqrt{z_2})=0\} .
$$
The holomorphic polynomial map $$F: (w_1,w_2) \mapsto (z_1,z_2) =  (w_1,w_2^2)$$ 
satisfies
$F^{-1}(\G^*) = \Pi_1 \cup \Pi_2$,  where $\Pi_1 = \{ {\rm Im\,}(w_1 + w_2) = 0 \}$ and 
$\Pi_2 = \{ {\rm Im\,}(w_1 - w_2) = 0 \}$. From the discussion in Section~\ref{s.BR} we know that a neighbourhood of the
origin can be ``filled" by holomorphic discs attached to a neighbourhood of the origin in 
$\Pi_1 \cup \Pi_2$. 

Given a continuous CR-function $f$ on $\G^*$, the function $g = f\circ F$ is CR on $\Pi_1 \cup \Pi_2$.
By the Baoudendi-Tr\`eves approximation theorem (see, e.g.,~\cite{B}) and the discussion above, the function 
$g$ extends holomorphically to some neighbourhood of the origin, call the extension $G(w)$. Then the 
composition $G\circ F^{-1}$ defines a $1$-to-$2$ holomorphic correspondence that extends the function~$f$. 
But in fact, $G\circ F^{-1}$ is single-valued, which can be seen as follows.  The set $\{z_2=0\}$ is 
the ramification locus of the correspondence $F^{-1}$. For a small $r>0$ the curve
$$
[0,2\pi] \ni \theta \mapsto (-\sqrt{r} \sin (\theta/2), r e^{i\theta})\in \cx^2
$$
is contained in $\G^*$ and is a generator of the fundamental group of $\cx^2\setminus\{z_2=0\}$.
This, however, contradicts the fact that the restriction of $G\circ F^{-1}$ to $\G^*$ equals $f$.

It follows that any continuous CR-function on $\G^*$ extends holomorphically to a full neighbourhood
of the origin in the ambient $\cx^2$. This implies that no compact neighbourhood of the origin in 
$\G^*$ is holomorphically convex, in particular, not rationally or polynomially convex.
\end{ex}

\bigskip

\noindent{\bf Acknowledgment.} We thank I.~Kossovskiy for very useful discussions.


\end{document}